\documentclass{amsart}
\usepackage{url}
\usepackage{graphicx, float}

\usepackage{mathtools}
\usepackage{hyperref}
\newcommand{\R}{\mathbb{R}}
\newcommand{\N}{\mathbb{N}}

\newtheorem{theorem}{Theorem}[section]

\theoremstyle{definition}

\newtheorem{conjecture}[theorem]{Conjecture}

\theoremstyle{remark}

\numberwithin{equation}{section}

\begin{document}

\title[AMS Article Template]{Global and Local Bounds on the fundamental ratio of Triangles and Quadrilaterals}


%    author one information
%    author two information
\author{Ryan Arbon}

\address{Department of Mathematics, University of California: Los Angeles, Los Angeles, CA 90095, USA}
\email{rarbon@math.ucla.edu}
\curraddr{}
%    \subjclass is required.

%    Abstract is required.
\begin{abstract}
We present a new, computer-assisted, proof that for all triangles in the plane, the equilateral triangle uniquely maximizes the ratio of the first two Dirichlet--Laplacian eigenvalues. This proves an independent proof the triangular Ashbaugh--Benguria--Payne--P\'olya--Weinberger conjecture first proved in \cite{Siudeja} and \cite{Arbon}. Inspired by \cite{Lu}, the primary method is to use a perturbative estimate to determine a local optimum, and to then use a continuity estimate for the fundamental ratio to perform a rigorous computational search of parameter space. We repeat a portion of this proof to show that the square is a strict local optimizer of the fundamental ratio among quadrilaterals in the plane.
\end{abstract}

\maketitle

\pagestyle{plain}

\bibliographystyle{amsplain}

\section{Introduction}

Let $D$ be an bounded open domain in the plane $\R^2$, and let $\{\lambda_i\}_{i \in \N}$ denote the sequence of Dirichlet--Laplacian eigenvalues of the Laplacian $\Delta$ on $D$, ordered such that $\lambda_1 \leq \lambda_2 \leq \lambda_3 \leq \hdots \longrightarrow \infty$. We shall refer to the $\lambda_i$ as the eigenvalues of $D$ throughout this paper.\\

Famously, the collection $\{\lambda_i\}$ does not uniquely determine $D$. However, there has been interest in finding optimizing domains $D$ for various functions of the $\lambda_i$, including the \textbf{fundamental gap}:

\begin{equation*}
    \chi(D) \coloneqq \lambda_2(D) - \lambda_1(D)
\end{equation*}

and the \textbf{fundamental ratio}:

\begin{equation}\label{spectralratio}
    \xi(D) \coloneqq \frac{\lambda_2}{\lambda_1}(D)
\end{equation}

Our focus in this paper will be proving bounds on \eqref{spectralratio} when $D$ is restricted to polygonal subclasses of $2^{\R^2}$. In 1955, Payne, P\'olya, and Weinberger proved in \cite{PPW} that among general $D \subset \R^2$, we have $\xi(D) \leq 3$. They furthermore conjectured that the unique optimizer of $\xi$ is the disk. Later, Thompson proved analogous bounds in dimension $d$ and conjectured that among domains $D \subset \R^d$, the optimizer is the $d$-dimensional ball \cite{thompson}.\\

The Payne--P\'olya--Weinberger (PPW) Conjecture, as it came to be known, was successfully proven by Ashbaugh and Benguria in 1991 and 1992 in \cite{Ashbaugh1991ProofOT, Ashbaugh1992}. Hence, one expects that domains $D$ which are somehow close to the disk will provide near optimal values of $\xi$. From another perspective, we expect domains with high amounts of symmetry to have large values of $\xi$. Restricting $D$ to lie in a fixed class of $k$-sided polygons leads to the polygonal Ashbaugh--Benguria--PPW conjecture, which was first stated by Antunes and Freitas in \cite{AntFrei}, who conducted a numerical study of the conjecture. The general case remains open, and we state as in \cite{open, henrot2017shape}:

\begin{conjecture}\label{bigPPW}
Let $P_k \subset \R^2$ denote an open $k$-gon and let $P_k^* \subset \R^d$ denote an open regular $k$-gon. Then for $k \geq 3$:

\begin{equation}\label{bigPPWeqn}
    \xi(P_k) \leq \xi(P_k^*)
\end{equation}
\end{conjecture}

For $k = 3$, this inequality was proved for acute triangles by Siudeja in 2010 in \cite{Siudeja}, and the remaining $k=3$ cases were proven by Arbon et al. in \cite{Arbon} in 2020.

The primary contribution of this paper is to provide a new proof of the $k = 3$ case of \ref{bigPPW}, with the additional property that the optimizer is explicitly shown to be unique:

\begin{theorem}\label{triPPW}
For any triangle $P_3 \subset \R^2$, we have

\begin{equation*}
    \xi(P_3) \leq \xi(P_3^*)
\end{equation*}

for $P_3^*$ an equilateral triangle. Furthermore, equality holds if and only if $P_3$ is equilateral.
\end{theorem}

To achieve this result, we first show that an appropriately normalized equilateral triangle is a strict local optimum using eigenvalue perturbation methods as in \cite{Courant} to establish a local bound.  Next, we use translation, rotation, and reflection invariance of the fundamental ratio to establish a bounded parameter space, which we compactify to prevent degeneration of the triangles to a line. We then establish a continuity estimate for $\xi$ in the moduli space of triangles, and we use this estimate to perform a computational search of the remaining parameter space at a finite number of points. The search is performed using finite element software, and we use bounds as proved by Liu in \cite{Liu} to perform validated numerical computation. This approach is inspired by Lu and Rowlett, who used similar ideas in \cite{Lu} to prove that the fundamental gap is uniquely minimized by the equilateral triangle among the moduli space of appropriately normalized triangles.\\

The final contribution of this paper is a re-appliation of the previously mentioned perturbative techniques to prove that the square is a strict local maximizer of $\xi$ among qadrilaterals (with the appropriate parametrization):

\begin{theorem}\label{QuadThm}
For $P_4$ sufficiently close to $P_4^*$, we have

\begin{equation*}
    \xi(P_4) \leq \xi(P_4^*)
\end{equation*}

for $P_4^*$ a square. Furthermore, equality holds if and only if $P_4$ is square.
\end{theorem}

Lastly, we discuss attempts to generalize the computational approach for triangles to the quadrilateral cases, together with the accompanying difficulties. The key problem is in managing degenerate cases; while after normalization, triangular regions can only degenerate in one way, quadrilaterals can degenerate in several distinct manners, including non-convex behavior. Restricting to only convex quadrilaterals may simplify the proof.

\section{Triangles}

We will now show that the equilateral triangle is the unique maximizer for (\ref{bigPPW}) in the moduli space of triangles, starting with the local case.

\subsection{Local Maximality for Triangles}\label{localmax}

Let $T$ be the equilateral triangle with vertices $(0,0) $, $(1,0)$, and $(\frac{1}{2},\frac{\sqrt{3}}{2})$. We will write the eigenvalues of $T$ as $\lambda_i$. The equilateral triangle is one of the few domains in $\R^2$ with explicitly known eigenvalues. In particular, we have

\begin{equation}
    \xi(T) = \frac{\lambda_2}{\lambda_1}(T) = \frac{7}{3}
\end{equation}

By translation and rotation invariance, a perturbation of $T$ is determined by the perturbation of the $(\frac{1}{2},\frac{\sqrt{3}}{2})$ coordinate. We may either write the perturbed coordinate as $(p,q)$ or as $(\frac{1}{2},\frac{\sqrt{3}}{2}) + t(a,b)$ where $t\geq 0$ and $(a,b)$ is a unit vector. We will use the second formulation in this section, and for a fixed choice of $(a,b)$, we will write $T(t)$ for the perturbed triangle with eigenvalues $\lambda_i(t)$. The associated linear transformation taking $T(0) \to T(t)$ is represented by the matrix

\begin{equation}\label{amat}
    A = \begin{bmatrix}1 & t \frac{2a}{\sqrt{3}}\\
    0 & 1 + t \frac{2b}{\sqrt{3}}
    \end{bmatrix}
\end{equation}

As noted in \cite{Lu}, we have that $T(t)$ is isomorphic to $T$ with the metric

\begin{equation}
    g = (dx)^2 + \frac{4ta}{\sqrt{x}}dxdy + \left(\left(1+\frac{2tb}{\sqrt{3}}\right)^2 + \frac{4t^2a^2}{3}\right)(dy)^2.
\end{equation}

The corresponding Laplacian, again in similar notation to that of \cite{Lu}, is

\begin{equation}
    \Delta(t) = \Delta_0 - t \tilde{L}_1 - t^2 \tilde{L}_2
\end{equation}

where $\Delta_0 = \partial_x^2 + \partial_y^2$ and

\begin{equation}
    \tilde{L}_1 = \frac{4\sqrt{3}}{(\sqrt{3}+2tb)^2}(a\partial_x \partial_y + b \partial_y^2), \; \; \; \; \tilde{L}_2 = \frac{4}{(\sqrt{3}+2tb)^2}(-a^2\partial_x \partial_y + b \partial_y^2)
\end{equation}

We now wish to compute the perturbations of $\lambda_2/\lambda_1$ corresponding to this perturbed Laplacian. We first make a note about the symmetry of the problem.\\

Given a perturbation of strength $t \geq 0$ in the direction $(a,b)$, we will denote the fundamental ratio by

\begin{equation}
    \xi_{(a,b)}(t) = \frac{\lambda_2(t)}{\lambda_1(t)}
\end{equation}

We have only defined this for $t \geq 0$, so we extend our definition for negative arguments $-t <0$ as $\xi_{(a,b)}(-t) = \xi_{(-a,-b)}(t)$.

Consider now a deformation of strength $t$ in the direction $(1,0)$. Because the fundamental ratio is invariant under rotations, we have that $\xi_{(0,1)}(t) = \xi_{(0,-1)}(t) = \xi_{(0,1)}(-t)$. This holds for any $t$, no matter how small, so $\xi_{(0,1)}$ is an even function of $t$ around $t=0$. Hence the equilateral triangle is a critical point for perturbations, positive or negative, in the $(0,1)$ direction. Thus, if the fundamental ratio were differentiable in this direction, we would expect its first derivative to be zero.

Similarly, we now consider a perturbation of strength $t$ in the direction $(\frac{1}{2},\frac{\sqrt{3}}{2})$. By flipping the triangle around the line $y = \frac{x}{\sqrt{3}}$ and then shrinking the resulting triangle, we again obtain
$\xi_{(\frac{1}{2},\frac{\sqrt{3}}{2})}(t) = \xi_{(\frac{1}{2},\frac{\sqrt{3}}{2})}(-t)$ for this perturbation. Here we have exploited the scale invariance of the fundamental ratio. Thus the fundamental ratio is an even function of $t$ along $(\frac{1}{2},\frac{\sqrt{3}}{2})$. Since these are two linear independent perturbations, we expect that $t = 0$ is a critical point for any linear perturbation.

Because the second eigenspace of the equilateral triangle has dimension two, $\lambda_2(t)$ is not differentiable at $0$. Hence $\xi$ is not differentiable. Nevertheless, we shall see that if we fix $(a,b)$ and use $t \geq 0$, we can define a perturbative expansion of $\xi_{(a,b)}(t)$. The first order term in this expansion will be shown to be negative, and we will bound the remaining terms. In line with the symmetry of the problem, we will see that the first order estimate is equally valid for a perturbation $t$ in the direction $(a,b)$ as it is for $t$ in the direction $(-a,-b)$.\\

Visually speaking, if we view the graph of $\xi$ over $(p,q)$ space, we would see that along lines passing through $(1/2, \sqrt{3}/2)$, $\xi$ is symmetric (at least to first order). Visually, $\xi$ has a corner similar to the absolute value function. We proceed now with the perturbation theory.\\

We will use $u_i(t)$ to denote the $i$-th eigenvalue of $\Delta(t)$. For simple eigenvalues, we may expand around $t = 0$ to obtain

\begin{equation}
    \begin{split}
        u_i(t) &= u_i + t v_i + \mathcal{O}(t^2)\\
        \lambda_i(t) &= \lambda_i + t \nu_i + \mathcal{O}(t^2)
    \end{split}
\end{equation}

The perturbation of degenerate eigenvalues is more delicate. However, we will still use the symbol $\nu_i$ to denote the first order perturbation of such values.  Using this notation, the first order perturbation of the fundamental ratio becomes

\begin{equation}\label{chain}
    \frac{\nu_2 \lambda_1 - \nu_1 \lambda_2}{\lambda_1^2}.
\end{equation}

This formulation is valid since we are restricting ourselves to positive perturbations in a single direction. With this fixed direction, the second eigenspace splits into two fixed eigenspaces. Hence we may take a limit as $t \to 0^+$ and obtain Eq. (\ref{chain}).

We begin by finding $\nu_1$ and $\nu_2$. Due to the degeneracy of the second eigenspace, we will simultaneously also find $\nu_3$.

Consider $\lambda_1$ and $u_1$. We have

\begin{equation}
    \Delta_0 u_1(t) - t\frac{4\sqrt{3}}{(\sqrt{3}+2tb)^2}(a\partial_x \partial_y + b \partial_y^2) u_1(t) -t^2 \frac{4}{(\sqrt{3}+2tb)^2}(-a^2\partial_x \partial_y + b \partial_y^2)u_1(t) = -\lambda_1(t)u_1(t)
\end{equation}

The $\frac{1}{(\sqrt{3}+2tb)^2}$ term could present complications. If we take $t \in [0,0.8]$, we may expand around $t=0$ for any $b$ to obtain

\begin{equation}
    \Delta(t) = \Delta_0 - t\frac{4\sqrt{3}}{3}(a\partial_x \partial_y + b \partial_y^2) -t^2\frac{4}{3}((-a^2-4ab)\partial_x \partial_y - 3b^2 \partial_y^2) + \mathcal{O}(t^3)
\end{equation}

We therefore define

\begin{equation}
    L_1 = \frac{4\sqrt{3}}{3}(a\partial_x \partial_y + b \partial_y^2),  \; \; \; \; L_2 = \frac{4}{3}((-a^2-4ab)\partial_x \partial_y - 3b^2 \partial_y^2).
\end{equation}

In general we will write $L_n$ for the $n$-th order expansion. This absolutely converges because of our restrictions on $t$. Then to compute $\nu_1$ we collect the first order terms of

\begin{equation}
    \left(\Delta_0-tL_1 - t^2L_2 -\sum_{n=3}^\infty t^n L_n\right)u_1(t) = -\lambda_1(t)u_1(t)
\end{equation}

to obtain

\begin{equation}
    \Delta_0 + \lambda_i v_1 = L_1 u_1 - \nu_1 u_1
\end{equation}

which, as in \cite{Courant}, enables us to find

\begin{equation}
    \nu_1 = \int_T L_1 u_1 u_1.
\end{equation}

The normalized form for $u_1$ is

\begin{equation}
    u_1(x,y) = \frac{2\sqrt{2}}{3^{3/4}}\sin\left(\frac{2 \pi \sqrt{3} y}{3}\right) \sin\left(\pi\left(x + \frac{\sqrt{3}y}{3}\right)\right)\sin\left(\pi\left(x - \frac{\sqrt{3}y}{3}\right)\right),
\end{equation}

and so we compute

\begin{equation}
    \nu_1 = \frac{4\sqrt{3}}{3} b \left(-\frac{8 \pi^2}{3}\right).
\end{equation}

This first order result should not come as a surprise. Note that a positive $b$ value (and positive $t$) stretches the triangle upwards. Because of domain monotonicity, we would expect such a deformation to decrease $\lambda_1$. Similarly, we would expect a deformation downward to increase $\lambda_1$.\\

We now turn to the case of $\lambda_2(t)$. If we were interested in finding the perturbations of the second eigenvector itself, we would need to first seek $u_2^* =\alpha u_2 + \beta u_3$, with $\alpha$ and $\beta$ chosen such that the perturbation is continuous. As a vector, $(\alpha,\beta)$ is an eigenvector of the matrix

\begin{equation}
    D = \begin{bmatrix}d_{22} & d_{23}\\ d_{32} & d_{33}\end{bmatrix}
\end{equation}

where 

\begin{equation}
    d_{ij} = \int_T L_1 u_i u_j.
\end{equation}

However, we do not need to explicity compute $\alpha$ and $\beta$, since we are only concerned with the perturbation of the eigenvalues, rather than perturbation of the eigenvectors. The eigenvalues of the matrix $D$ are in fact the desired perturbations, $\nu_2$ and $\nu_3$.

Computing $D$ explicitly gives

\begin{equation}
    D = \frac{4\sqrt{3}}{3} \begin{bmatrix}a \frac{6561 \sqrt{3}}{800}  - b\frac{59049+44800 \pi^2}{7200}& -a \frac{6561}{800} - b \frac{6561\sqrt{3}}{800}\\ -a \frac{6561}{800} - b \frac{6561\sqrt{3}}{800} & -a \frac{6561 \sqrt{3}}{800}  + b\frac{59049-44800 \pi^2}{7200} \end{bmatrix},
\end{equation}

and a simple calculation using the fact that $a^2+b^2 = 1$ shows that the eigenvalues are $\frac{4\sqrt{3}}{3} \frac{ -59049 - 22400 b \pi^2}{3600}$ and $\frac{4\sqrt{3}}{3} \frac{59049 - 22400 b \pi^2}{3600}$. One of these eigenvalues is smaller than the other for any admissible $b$, and hence we denote

\begin{equation}
    \nu_2 = \frac{4\sqrt{3}}{3} \frac{ -59049 - 22400 b \pi^2}{3600}, \; \; \nu_3 = \frac{4\sqrt{3}}{3} \frac{59049 - 22400 b \pi^2}{3600}
\end{equation}

Note that we have this distinction for any $b$. Hence in our given direction, the second eigenspace splits in the perturbation, with the new second eigenspace being $\lambda_2(t) = \lambda_2 + t \nu_2 + \mathcal{O}(t^2)$. Then applying the form of the chain rule in this specific direction - where the splitting of the eigenspace is fixed, we find that

\begin{equation}
    \frac{\nu_2 \lambda_1 - \nu_1 \lambda_2}{\lambda_1^2} = - \frac{4\sqrt{3}}{3}\frac{19683}{6400 \pi^2} \leq -0.719632
\end{equation}

Meaning that

\begin{equation}
    \xi_{(a,b)}(t) = \frac{7}{3} -\frac{6561 \sqrt{3}}{1600 \pi^2} t + \mathcal{O}(t^2).
\end{equation}

If we can estimate the $\mathcal{O}(t^2)$ term as $Ct^2$, we will establish that $\lambda_2/\lambda_1$ is locally maximized by the equilateral triangle and we may estimate the region where the local maximization holds.

We can obtain a very rough estimate of $\lambda_1(t)$ by noting that for any $t \in [0,1/2]$, we may enclose $T(t)$ in a rectangle with side lengths $1$ and $\frac{1+\sqrt{3}}{2} $, and so by domain monotonicity, we have $\lambda_1(t) \geq \pi^2 (1 + \frac{4}{(1+\sqrt{3})^2}) \geq 15.1586$ for all such $t$ and any direction. We obtain a similar estimate of $\lambda_2(t)$ by noting that we may enclose a rectangle of side lengths $1/2$ and $1/10$ inside any $T(t)$ for $t\in[0,1/2]$, and so we have $\lambda_2(t) \leq 116 \pi^2 \leq 1144.88$.

Thus we have

\begin{equation}
    \xi_{(a,b)}(t) = \frac{7}{3} -\frac{6561 \sqrt{3}}{1600 \pi^2} t + \mathcal{O}(t^2) \leq \frac{1144.88}{15.1586} \leq 75.53
\end{equation}

and so the second order term satisfies $\mathcal{O}(t^2) \leq 295 t^2$ for $t \in [0,1/2]$, meaning that

\begin{equation}\label{triangleloc}
    \xi_{(a,b)}(t) \leq \frac{7}{3} -\frac{6561 \sqrt{3}}{1600 \pi^2} t + 295 t^2,
\end{equation}

and so $\xi_{(a,b)}(t) \leq \frac{\lambda_2}{\lambda_1}$ for $t \in [0,0.0022]$, and the equality is only realized for $t = 0$. This holds for any positive perturbation in any choice of direction, and hence the equilateral triangle is locally maximal. Note that inherent in \eqref{triangleloc}, we see that $T$ is a strict local maximizer.

\subsection{Extremal Cases and the Moduli Space of Triangle}

By the scale invariance of $\xi$, for we may without loss of generality assume that any triangle $T$ to have a longest side of length $1$ and to be aligned with the $x$-axis in the Euclidean plane. Hence we need only consider triangles with vertices $(0,0)$, $(1,0)$, and $(p,q)$, where $p$ and $q$ satisfy $p^2 + q^2 \leq 1$. By reflection invariance, we may also impose $p \geq 1/2$, $q > 0$. Thus after removing the symmetries, the only equilateral triangle is one considered in \ref{localmax}. Hence the equilateral triangle is a strict local maximizer in this parameter space.

In the following sections, will undertake a numerical search of the parameter space. However the numerical method cannot accurately resolve degenerating triangle. Under our given conditions, there is only one way in which a sequence of triangles can degenerate, which is for $q$ to go toward zero.

A rough estimate for $\xi$ in this scenario was given in \cite{Arbon}, where it was shown that

\begin{equation}\label{bound1}
    \xi(p,q) \leq \frac{(1 + \sqrt[3]{4q^2})^3}{(q+1)^2}
\end{equation}

and hence $\xi < \frac{7}{3}$ for $q \leq 0.156$. Together with the local minimality result, we have carved out two regions in $(p,q)$ space where we know Theorem \ref{triPPW} holds. The remaining set of $(p,q)$ points form a compact, connected subset of $\R^d$. In the following section, we will use numerical techniques to rigorously show that the inequality holds in this region as well.

\subsection{Continuity Estimate}\label{maximality}

While we used the matrix $A$ given by Eq. (\ref{amat}) to transform between the equilateral triangle and a perturbation of that triangle, a more general form of $A$ is described in \cite{Lu}:

\begin{equation}
    A = \begin{bmatrix}1 & \frac{t a}{q}\\
    0 & 1 + \frac{t b}{q}\end{bmatrix}
\end{equation}

This matrix gives a linear transformation from a triangle $T$ with vertex $(p,q)$ to a perturbed triangle $T'$ with vertex $(p,q) + t(a,b) = (p',q')$. Associated to this matrix $A$ is a change of metric $g$. The inverse metric $g^{-1}$ has eigenvalues $\gamma_-$ and $\gamma_+$ which gives the following bound as in \cite{Lu}

\begin{equation}
    |\lambda_i' - \lambda_i| \leq (\gamma_+ - \gamma_-)\lambda_i
\end{equation}

where $\lambda_i$ is an eigenvalue of $T$ and $\lambda_i'$ is an eigenvalue of $T'$. This follows from the fact that $\gamma_- \lambda_i \leq \lambda_i' \leq \gamma_+ \lambda_i$.

Plugging in the appropriate values and simplifying gives the bound on $\xi$:

\begin{equation}
\begin{split}
        |\xi' - \xi| &= \frac{\lambda_1|\lambda_2'-\lambda_2| + \lambda_2|\lambda_1 - \lambda_1'|}{\lambda_1 \lambda_1'}\\
        &\leq \frac{\lambda_1 + \lambda_2}{\lambda_1 \lambda_1'}(\gamma_+ - \gamma_-)\\
        &= \left(\frac{1+ \xi}{\lambda_1'} \right)(\gamma_+ - \gamma_-)\\
        &= \left(\frac{1+ \xi}{\lambda_1'} \right) t \frac{\sqrt{4q(q+tb)+t^2}}{(q+tb)^2}
\end{split}
\end{equation}

Now if we have only calculated $\xi$, we do not know $\lambda_1'$. However, since our triangles lie in the rectangle with sides lengths $1$ and $\sqrt{3}/2$, by domain monotonicity we can estimate $\lambda_1' \geq 23$, and so

\begin{equation}\label{contest1}
\begin{split}
        |\xi' - \xi| \leq \left(1 + \xi \right) t \frac{\sqrt{4q(q+bt)+t^2}}{23 (q+tb)^2} \leq \left(1 + \xi \right) t \frac{\sqrt{4q(q+t)+t^2}}{23 q^2} 
\end{split}
\end{equation}

if we restrict ourselves to $t \geq 0$ and $b \geq 0$. So now if we calculate numerically $\xi(p,q)$ at some fixed $(p,q)$ where $\xi(p,q)<7/3$, we can determine a region $U$ in parameter space, containing $(p,q)$ such that $\xi(p',q') < 7/3$ for $(p',q')\in U$.

Suppose that $\xi(p,q) < 7/3$. Certainly $\xi(p',q') < 7/3$, for $p' = p+ta$, $q'=q+tb$ so long as

\begin{equation}
    \begin{split}
        \xi(p,q) + |\xi(p,q) - \xi(p',q')| < \frac{7}{3}
    \end{split}
\end{equation}

Then by using the continuity estimate Eq. (\ref{contest1}), we find

\begin{equation}\label{contest2}
    \begin{split}
        \xi(p,q) + |\xi(p,q) - \xi(p',q')| & \leq \xi(p,q) + (1+\xi(p,q)) \frac{\sqrt{4q(q+t)+t^2}}{23 q^2}
    \end{split}
\end{equation}

which holds so long at $b \geq 0$ and $t \geq 0$. If we choose a $t > 0$ such that the right hand side of Eq. (\ref{contest2}) is less than $7/3$, then we have determined a half circle of radius $t$ in parameter space, centered at the original $(p,q)$ satisfying the desired property. For the sake of numerical simplicity, we note that if we take $t^* = \frac{\sqrt{2}}{2}$ then we can determine a square of side-length $t^*$ whose bottom left corner is $(p,q)$ also satisfying this property.

Letting the right hand side of \eqref{contest2} be denoted as $f(p,q)$, we have

\begin{equation}
    f_t(p,q) = \xi(p,q) + \left(1+\xi(p,q)\right)t \frac{\sqrt{4q(q+t)+t^2}}{23 q^2}
\end{equation}

We wish to find values of $t$ such that $f_t(p,q) < 7/3$. Noting that $a,b,t,q,p > 0$, we see that equivalently we need to find $t$ where the following expression is equal to zero:

\begin{equation}\label{maximalt}
    \frac{(1 + \xi)}{23 q^2} t^2 + \frac{2(1+\xi)}{23 q}t + \xi - \frac{7}{3}
\end{equation}

If we compute the positive zero of Eq. (\ref{maximalt}), we find a maximal $t$ such that $\xi(p,q) + |\xi(p,q) - \xi(p',q')|\leq 7/3$. Then for such a $t$, if we take $t^* = \frac{\sqrt{2}}{2}t - \epsilon$ for some small $\epsilon >0$, the desired square property will hold, now with strict inequality.\\

\subsection{Validated Numerics}\label{valid}

Until this point, we have implicitly assumed that we are capable of exactly evaluating $\xi(T)$ for any triangle in our parameter space, and then using exact arithmetic to determine the $t^*$ associated to $T$. Unfortunately, this is not the case. However, we can bound $\xi(T)$ from above and then use validated numerics to obtain a comparable result.

In particular, we use a lower bound on $\lambda_1(T)$ and an upper bound on $\lambda_2(T)$ coming from a finite element method. In \cite{Liu}, Liu proved the following:

\begin{theorem}
    Let $D$ be a polygonal domain triangulated by a mesh of size at most $h$. Let $\lambda_k$ be the $k$-th eigenvalue of $D$ and let $\lambda_{k,h}$ be the $k$-th eigen value of $D$ as computed by the Crouzeix--Raviert method. Then
    
    \begin{equation}\label{l1bnd}
        \frac{\lambda_{k,h}}{1 + (0.1893 h)^2 \lambda_{k,h}} \leq \lambda_k
    \end{equation}
\end{theorem}

In our approach, we shall therefore use the Crouzeix--Raviert method to compute $\lambda_{1,h}$. This then gives the lower bound on $\lambda_1$ from \eqref{l1bnd}.

In the process of the finite element method, we will produce two orthonormal approximate eigenfunctions $u_1$, $u_2$ corresponding to $\lambda_1$ and $\lambda_2$, respectively. As is well known, the Rayleigh--Poincar\'e Principle gives a variational characterization of the sum of eigenvalues of a domain $D$:

\begin{equation}
    \lambda_1 + \hdots + \lambda_n = \min \{ \int_D |\nabla v_1|^2 dx + \hdots  + \int_D |\nabla v_n|^2 dx : v_j \in H_0^1(D), \langle v_j, v_k\rangle = \delta_{jk}\}
\end{equation}

We can therefore use the computed approximate eigenfunctions to give an upper bound on $\lambda_1 + \lambda_2$ as $\lambda_1 + \lambda_2 \leq \int_T |\nabla u_1|^2 dx + \int_T |\nabla u_2|^2 dx$. Then as $\lambda_2 = \lambda_1 + \lambda_2 - \lambda_1$, we have

\begin{equation}\label{l2bnd}
    \lambda_2 \leq \int_T |\nabla u_1|^2 dx + \int_T |\nabla u_2|^2 dx + \frac{\lambda_{1,h}}{1 + (0.1893 h)^2 \lambda_{1,h}}
\end{equation}

Combining the bounds in \eqref{l1bnd} and \eqref{l2bnd} gives an upper bound $\xi_h(T)$ for the true fundamental ratio of a triangle $\xi(T)$. Substituting $\xi_h(T)$ into the right side of \eqref{contest2} gives us

\begin{equation}
    \xi(p,q) + |\xi(p,q) - \xi(p',q')| \leq \xi_h(p,q) + (1+\xi_h(p,q)) \frac{\sqrt{4q(q+t)+t^2}}{23 q^2}
\end{equation}

Proceeding as before, we find that $t^*$ is now given by the positive zero of

\begin{equation}\label{ttrue}
    \frac{(1 + \xi_h)}{23 q^2} t^2 + \frac{2(1+\xi_h)}{23 q}t + \xi_h - \frac{7}{3}
\end{equation}

multiplied by $\sqrt{2}/2$.\\

Lastly, we mention that due to the imprecise nature of machine arithmetic, it is necessary at each step in the proceeding algorithm to account for the accumulation of machine and other errors. For instance, in the algorithm, we use $0.9*t^*$ rather than $t^*$, both in order to ensure strict inequality and to account for errors introduced by the quadratic formula. In most cases, we add or subtract terms of order $\epsilon = 1e-9$ in the direction that weakens the bound. While most steps of the algorithm are more accurate than this, the error is dominated by the computation of $\lambda_{1,h}$ (and $\lambda_{2,h})$ with an error tolerance of $1e-10$. Hence out of an abundance of caution, we use errors of order $\epsilon$ throughout the algorithm

\subsection{Algorithm}

By the work in the previous sections, we need only concern ourselves with calculating $\xi$ for triangles with vertices $(0,0)$, $(1,0)$, and $(p,q)$ such that

\begin{itemize}\label{conditions}
    \item $p^2 + q^2 \leq 1$,
    \item $0.5 \leq p \leq 1$,
    \item $0.156 \leq q \leq 1$ and
    \item $\sqrt{(p-1/2)^2 + (q-\sqrt{3}/2)^2} > 0.0022$.
\end{itemize}

The calculations in the algorithm are done by the FreeFEM++ software \cite{FEM}. The idea of the algorithm is to choose a point $(p_i,q_j)$ in the parameter space such that $\xi(p_i,q_i)<7/3$. Then by Subsections \ref{maximality} \ref{valid}, we find a square of side length $t^*$ such that the inequality holds (strictly) within it. Our next point will be another corner of the square. In terms of specifics, we begin with the point $(p_{0,0},q_0) = (0.5-\epsilon,0.156-\epsilon)$. We then proceed by advancing $p_0$ to $p_{1,0} = p_0 + t^*_0$. We proceed until $(p_{i,0} + t^*_i,q_j)$ would violate the above conditions on the parameter space. We then return to $p_{0,1} = p_{0,0}$ and advance $q_0$ to $q_1 + t^{**}_0$, where $t^{**}_0$ is the minimum of the $t^*_{i,0}$. This process then repeats. Algorithmically speaking, the steps are as follows:\\

0. Define $p_{0,j} = 0.5$ for all $j$ and $q_0 = 0.156$.\\

1. At step $(i,j)$, we compute $\lambda_{1,h}(p_{i,j},q_j)$ and approximate eigenfunctions $u_1(p_{i,j},q_j)$ and $u_2(p_{i,j},q_j)$ using the Crouzeix--Raviert method. After ensuring that $u_1$ and $u_2$ are orthonormal up to order $1e-16$, we compute $\xi_h(p_{i,j},q_j)$ and check that $\xi_h < 7/3$. At each step in this process, we add errors of size $\pm \epsilon$ acorrding to the direction which weakens the bound.\\

2. We calculate $t^*_{i,j}$ by finding the positive root of Eq. (\ref{ttrue}) and multiplying this root by $0.9*\sqrt{2}/2$. Then by our continuity estimate, any triangle defined by a vertex in the square of side length $t^*_{i,j}$ with lower left corner $(p_{i,j},q_j)$ strictly satisfies Theorem (\ref{triPPW}).\\

3. We now advance in $i$. We propose the point $p_{i+1,j} = p_{i,j} + t^*_{i,j}$. With the proposed point, we check $p_{i+1,j}^2 + q_j^2 \leq 1$ and $0.5 \leq p_{i+1,j} \leq 1$. We update $t^{**}_{i,j} = \min\left(t^{**}_{j},t^*_{i,j}\right)$. If both of the inequalites are satisfied, we accept $p_{i+1,j}$ and return to step 1. If either inequality is not satisfied, we move to step 4. Note that since we are advancing in $i$, we know that $p_{i,j}$ satisfies $\sqrt{(p_{i,j}-1/2)^2 + (q_j-\sqrt{3}/2)^2} > 0.0022$, so $p_{i+1,j}$ will as well.\\

4. We now attempt to advance in $j$. We propose the point $q_{j+1} = q_{j} + t^{**}_j$. We check $(p_{0,j+1},q_{j+1})$ to see if it lies in our parameter space. If the only inequality violated is $\sqrt{(p_{0,j+1}-1/2)^2 + (q_{j+1}-\sqrt{3}/2)^2}$ being $\leq 0.0022$, we re-define $p_{0,j+1} = \sqrt{0.0022^2 -(q_{j+1}-\sqrt{3}/2)^2} + 1/2$ and return to step 1. If any other inequality is violated, we are done and the program terminates. If none of the inequalities are violated, then we accept the point and return to step 1.\\

\subsection{Results}

\begin{figure}
    \centering
    \includegraphics[scale=0.6]{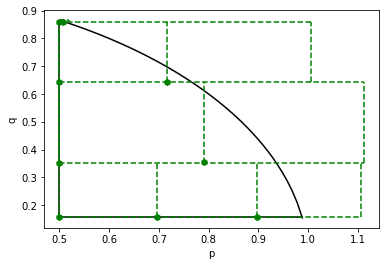}
    \caption{Representation of admissible $(p,q)$ parameter space, outlined in black, together with output of the triangle-checking algorithm, in green.}\label{triangleoutput}
\end{figure}

Running this algorithm completes the proof in the case of triangles. The FreeFEM++ code implementing the algorithm is all is available on GitHub\footnote{\href{https://github.com/rarbon/ComputationalPPWProof}{\texttt{https://github.com/rarbon/ComputationalPPWProof}}}. We can gain a visual sense of the output of the algorithm by examining Figure (\ref{triangleoutput}). Outlined in black, we see the region of $(p,q)$ space that we need to check Theorem (\ref{triPPW}) on. The green dots represent the points in the parameter space checked by the algorithm. Each green dot $(p_{i,j},q_j)$ lies at the corner of a rectangle of size $t^*_{i,j} \times t^{**}_{j}$ (visualized in dashed green), within which the inequality holds.

\begin{figure}
    \centering
    \includegraphics[scale=0.6]{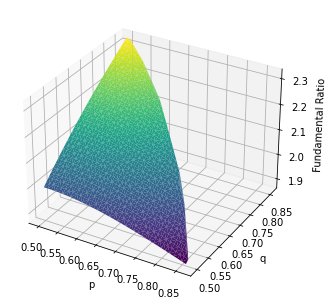}
    \caption{Fine grid of points in $(p,q)$ parameter space showing $\xi(p,q)$ as computed by FreeFEM++. Color scaling for visibility purposes.}\label{triangleratio}
\end{figure}

We see that relatively few points are needed to check this inequality, owing to both the strength of the estimate in (\ref{contest2}) and the low-dimensionality of the problem. For visual reference, we provide Figure (\ref{triangleratio}), where one sees that $\xi(p,q)$ decreases quickly away from $(1/2,\sqrt{3}/2)$, further providing intuition as to why so few evaluation points were needed. In fact, Figure (\ref{triangleratio}) appears to imply the $\xi$ is a concave function on the parameter space:

\begin{conjecture}\label{conv}
The fundamental ratio $\xi$ is concave on the moduli space of triangles once the symmetries of the problem have been remove.
\end{conjecture}

Coupled with the local maximization portion of the proof, a proof of Conjecture \ref{conv} would immediately imply Theorem \ref{triPPW} without the need for a computational parameter search.

\begin{figure}
    \centering
    \includegraphics[scale=0.6]{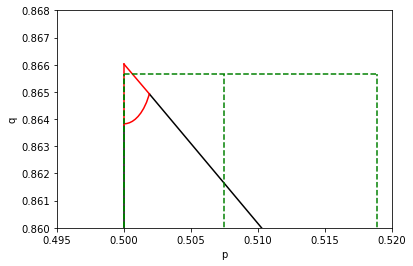}
    \caption{Zoom in of admissible $(p,q)$ parameter space, outlined in black, together with output of the triangle-checking algorithm, in green near $(1/2,\sqrt(3)/2)$. The excised region is outlined in red.}\label{trianglezoom}
\end{figure}

To see why we needed to excise a region around $(1/2,\sqrt{3}/2)$, see Figure (\ref{trianglezoom}), which is a zoomed in version of Figure (\ref{triangleoutput}), which the additional feature that the excised neighborhood of $(1/2,\sqrt{3}/2)$ is outlined in red.

Attempting to advance the algorithm by checking a point on the edge of the checkered rectangle, would result in a $t_{i,j}$ value of approximately $0.0009$. Further attempts to advance the algorithm would give smaller and smaller $t$ values, getting closer and closer to the critical point but never reaching it. Eventually, the limits of machine precision would either prevent the creation of more $t$ values or would make it impossible to trust the computer calculated $\xi$. Thus excising the neighborhood of $(1/2,\sqrt{3}/2)$ was critical in the algorithm's success.

\section{Quadrilaterals}

We now turn to the case of quadrilaterals. Alongside special triangles, rectangles are the only polygonal domains with known explicit forms for the eigenvalues. Hence we may hope that the problem may be amenable to more elementary techniques.

In the quadrilateral case, we may write Eq. (\ref{bigPPWeqn}) as

\begin{equation}\label{ppwquad}
    \xi(Q) = \frac{\lambda_2}{\lambda_1}(Q) \leq \frac{5 \pi^2}{2 \pi^2} = \frac{5}{2}
\end{equation}

where $Q$ ranges over all quadrilaterals. Our strategy in the proof of Theorem \ref{QuadThm} is very similar to the triangular proof. There is potential for a full global proof of Conjecture \ref{bigPPW} in the style of the proof of Theorem \ref{triPPW}. The primary difficulties lie in the increased dimensionality of the problem, as we will see that parameter space will require four coordinates rather than two. This increases both computation time as well as the number of extremal cases.

\subsection{Local Maximality for Quadrilaterals}\label{localproof}

We proceed in a similar fashion as in the case of triangles. Because linear transformations can only transport parallelograms to parallelograms, we rely on perspective mappings, which complicate the form of the Laplacian associated with a deformation. Nevertheless, the basic structure of the proof is the same.

Let $Q$ be the unit square with vertices at $(0,0)$, $(1,0)$, $(1,1)$, and $(0,1)$. Due to the scale invariance of the fundamental ratio, we may dilate and rotate any deformation, so we may assume without loss of generality that any local deformation of $Q$ will have coordinate vertices at $(0,0)$, $(1,0)$ thereby fixing the bottom side, with the remaining vertices being of the form $(0,1) + t(a,b)$, $(1,1) + s (c,d)$, where $(a,b)$ and $(c,d)$ are unit vectors, while $s$ and $t$ are non-negative scalars. In fact, we will see that we can assume that the side with endpoints $(0,0)$ and $(1,0)$ is the longest side. However, we will not explicitly require this at this point in time. Negative values of $s$ and $t$ would simply involve switching the direction of the respective perturbation(s).

Given a choice of $(a,b)$ and $(c,d)$, let $Q(t,s)$ denote the perturbed quadrilateral in these directions with respective strengths $t$ and $s$. The associated fundamental ratio could be denoted $\xi_{(a,b),(c,d)}(t,s)$, but for ease of notation we suppress the subscript where possible and write $\xi(t,s)$. We extend $\xi$ to negative arguments by defining $\xi_{(a,b),(c,d)}(-t,s) = \xi_{(-a,-b),(c,d)}(t,s)$ and $\xi_{(a,b),(c,d)}(t,-s) = \xi_{(a,b),(-c,-d)}(t,-s)$. %We will see that $(0,0)$ is a critical point of the fundamental ratio, and that in fact in any fixed direction, $\xi(t,s)$ is an even function in each argument up to first order. Ultimately, our work will be independent of the choice of $(a,b)$ and $(c,d)$.

Due to the symmetries of the quadrilateral, we notice that $\xi(t,s)$ possesses certain "evenness" properties, which would imply that in $(p_1,q_1,p_2,q_2)$ space, where the perturbable corners have coordinates $(p_1,q_1)$ and $(p_2,q_2)$, that $\nabla\xi (0,0,0,0) = 0$ (or in any fixed direction $\partial_i \xi(0,0) = 0$ in $(t,s)$ space), provided that $\xi$ was differentiable. Again we note that $\xi$ is not differentiable. Instead, $\xi$ has a corner at the regular quadrilateral. Geometrically, we will see that the first order perturbation of $\xi$ along a fixed direction pair $(a,b)$, $(c,d)$ is even, in the sense that for $t \geq 0$, $\xi(t,s) = \xi(-t,s)$ and for $s\geq 0$, $\xi(t,s) = \xi(t,-s)$ up to first order. We now proceed with the proof.\\

For a fixed $(a,b)$ and $(c,d)$, the perspective transformation taking $Q(t,s)$ to $Q$ is given by

\begin{equation}
    (x,y) \mapsto \frac{(a_0 x, a_1 y)}{(a_0+a_1-1)+(1-a_1)x+(1-a_0)y}
\end{equation}

where

\begin{equation}
    a_0 = 1 + s c - t \frac{1+s d}{1+t b}
\end{equation}

and

\begin{equation}
    a_1 = \frac{1+d t}{1+ tb}
\end{equation}

For more information on perspective mappings between convex quadrilaterals, see \cite{perspective}. This map is an immersion, and so we may pull back the Euclidean metric under it to obtain a new metric $g$. Written in terms of $a_0$ and $a_1$, we can represent $g$ in matrix form, which we write in components:

\begin{equation}
    g_{1,1} = \frac{a_0^2 a_1^2(-1+a_0+a_1)^2((2-2a_0+a_0^2)y^2 - 2a_1y(1-a_0+y)+a_1^2(1+y^2))}{((-1+a_1)a_1 x + a_0(a_1-y)+a_0^2 y)^4}
\end{equation}

\begin{equation}
    g_{1,2}=g_{2,1} = \frac{a_0^2 a_1^2(-1+a_0+a_1)^2(-a_1^2 x y + (a_0-2x+2a_0x-a_0^2x)y + a_1(x-a_0x-a_0y+2 x y))}{((-1+a_1)a_1 x + a_0(a_1-y)+a_0^2 y)^4}
\end{equation}

\begin{equation}
g_{2,2}=\frac{a_0^2 a_1^2(-1+a_0+a_1)^2((2-2a_1 + a_1^2)x^2 - 2 a_0 x(1-a_1+y)+a_0^1(1+x^2))}{((-1+a_1)a_1 x + a_0(a_1-y)+a_0^2 y)^4}
\end{equation}

The Laplacian associated with this metric is $\Delta = \frac{1}{\sqrt{\det{g}}}\sum_{i,j=1}^2 \partial_i g^{ij} \sqrt{\det{g}} \partial_j$, where $g$ with indices up represents the inverse of $g$. Computing this explicitly gives a lengthy expression involving several hundred different terms, along with a common denominator that is also a polynomial of high order. However, we need only consider the lowest order terms. After Taylor expanding the denominator, we find that for sufficiently small $t$ and $s$,

\begin{equation}
    \begin{split}
        \Delta = \Delta_0 - t L_t - s L_s + \mathcal{O}(2)
    \end{split}
\end{equation}

where

\begin{equation}
    L_t = 2\left(b \partial_x + a \partial_y + (ax+by)\partial_x\partial_y + (ay - (1-2x)b)\partial_x^2 + (bx-(1-2y)a)\partial_y^2\right),
\end{equation}

\begin{equation}
    L_s = -2\left(d \partial_x + c \partial_y + (cx+dy)\partial_x\partial_y + (cy - (1-2x)d)\partial_x^2 + (cx-(1-2y)d)\partial_y^2\right),
\end{equation}

\begin{equation}
    \Delta_0 = \partial_x^2 + \partial_y^2,
\end{equation}

and $\mathcal{O}(2)$ represents all terms of order $2$ and higher. A careful analysis of the necessary Taylor expansions reveals that $(t,s) \in [0,0.4]\times[0,0.4]$ is sufficiently small. Note that at first order, the perturbations are decoupled, and so the computation of the first order perturbations in $t$ and $s$ follow the same manner as in the case of triangles. We use $\nu_j^t$ and $\nu_j^s$ to denote the first order perturbation of the $j$-th eigenvalue corresponding to $t$ and $s$, respectively. Because the first eigenspace is simple, we find that

\begin{equation}
    \nu_1^t = \int L_t u_1 u_1 = - \pi^2(a+b)
\end{equation}

and

\begin{equation}
    \nu_1^s = \int L_s u_1 u_1 = \pi^2(c + d)
\end{equation}

Because the second eigenspace is not simple, we must proceed with more care. Again we form the matrices $D_t$ and $D_s$, with

\begin{equation}
    D_t = \begin{bmatrix}d_{22}^t & d_{23}^t\\ d_{32}^t & d_{33}^t\end{bmatrix} \; \; \; \; D_s = \begin{bmatrix}d_{22}^s & d_{23}^s\\ d_{32}^s & d_{33}^s\end{bmatrix}
\end{equation}

where $d_{ij}^t = \int L_t u_i u_j$ and $d_{ij}^s = \int L_s u_i u_j$. The eigenvalues of these matrices give the respective perturbations, and hence we compute that

\begin{equation}
\begin{split}
        \nu_2^t &= -\frac{5}{2}(a+b)\pi^2 - \frac{1}{18}\sqrt{16384(a+b)^2+729(a-b)^2 \pi^4}\\
        \nu_2^s &= \frac{5}{2}(c+d)\pi^2 - \frac{1}{18}\sqrt{16384(c+d)^2+729(c-d)^2 \pi^4}
\end{split}
\end{equation}

Now the first order perturbations of $\lambda_2/\lambda_1(t,s)$ in $t$ or $s$ will follow the same form as the quotient rule for differentiation, because we are focusing on a positive perturbation in a fixed direction. Thus we immediately have

\begin{equation}
    \frac{\lambda_2}{\lambda_1}(t,s) = \frac{5}{2} - t \frac{\sqrt{16384(a+b)^2+729(a-b)^2\pi^4}}{36\pi^2}- s \frac{\sqrt{16384(c+d)^2+729(c-d)^2\pi^4}}{36\pi^2} + \mathcal{O}(2).
\end{equation}

From this form, we determine that for any choice of directions, the fundamental ratio will decrease under a small enough perturbation. Furthermore, we can see from this formula that given a direction, the fundamental ratio is symmetric in the opposite direction. Hence this is a critical point and indeed a local maximum.

We now perform a crude estimate for $\lambda_2(t,s)$ and $\lambda_1(t,s)$ with $(t,s) \in [0,0.25] \times [0,0.25]$. We may enclose $Q(t,s)$ within a square of side length $3/2$ and we may enclose a square of side length $1/2$ within $Q(s,t)$. Hence $\lambda_2(t,s) \leq 20\pi^2$ and $\lambda_1(t,s) \geq 8 \pi^2/9$, implying that $\lambda_2/\lambda_1(t,s) \leq 45/2$.

Since $(a,b)$ and $(c,d)$ are unit vectors, the order one perturbations always have strength in between $0.509475$ and $1.060661$. Thus we seek constants $C_t$, $C_s$, $C_{ts}$ such that

\begin{equation}
    -1.060661(t+s) + C_tt^2 +C_ss^2 + C_{ts}ts \geq 20
\end{equation}

along the lines $t=0$, $s=1/4$ and $t=1/4$, $s=0$. Under these conditions, we may take $C_t = C_s = 336.972$. In fact, it is straightforward to see that we may also take $C_{ts} = 0$. Thus we have established that for any choice of directions

\begin{equation}
    \frac{\lambda_2}{\lambda_1}(t,s) \leq \frac{5}{2} -0.509475(t+s) + 336.972(s^2+t^2)
\end{equation}

and hence the fundamental ratio is locally maximized at $(0,0)$ and is the unique maximizer in the region $(t,s) \in [0,0.0015]\times[0,0.0015]$ for any choice of directions.

\subsection{Potential Parametrizations}\label{extra}

The natural parametrization following Subsection \ref{localproof} is to use scale invariance, rotation invariance, and translation invariance to ensure that every quadrilateral under consideration have longest side length 1 with vertices at $(0,0)$ and $(1,0)$. The two coordinates which may vary are labelled as $(p_1,q_1)$ and $(p_2,q_2)$.

While this parametrization does ensure that the parameter space is bounded, there are many possible forms of degeneration, the majority of which involve non-convex quadrilaterals. While initial numerical result seem to indicate that $\xi$ is well-behaved in the degenerating regimes, analytic bounds have not been established. If we restrict our attention to convex quadrilaterals, we may simplify the situation considerably.\\

For convex quadrilaterals, we will shift our perspective slightly. We will take as our representation of the regular quadrilateral the square with vertices at $(0,0)$, $(1/2, 1/2),$ $(1,0)$, $(1/2,-1/2)$. The two vertices which we may vary are $(1/2,1/2)$ and $(1/2,-1/2)$, which we label using $(p_1,q_1)$ and $(p_2,q_2)$ If we follow a similar proof as in Subsection \ref{localproof}, we may establish that this quadrilateral is a local maximizer with the region of local maximization scaled by $\sqrt{2}/2$.\\

Then for quadrilaterals under consideration, we enforce that the longest diagonal is the line segment joining $(0,0)$ and $(1,0)$. Under this regime, it is easy to examine the behavior of convex quadrilaterals. We have $q_2 \leq 0$ and $q_1 \geq 0$. One can also assume using reflection invariance that $|q_2| \leq q_1$ and $p_2 \leq 1/2$. These simplifaction provide for a regime where a parameter space search would be easy to organize.

The main obstacles to completing such a parameter space search at this time and hence establishing the polygonal Ashbaugh--Benguria--PPW inequality for quadrilaterals are continuity estimates and bounding the extremal cases. Such results seem possible by combining the correct bounds, as the literature contains many results, including asymptotics, for convex quadrilaterals. See for instance \cite{Freitas}.

\section*{Acknowledgement}

This work originated as a senior thesis project at Princeton University, with Peter Sarnak as primary advisor. The idea of the proof -- proving a local estimate then performing a computational search of a finite dimensional parameter space -- is due to him. Javier G\'omez-Serrano served as secondary advisor, and his comments contributed to making the computational portion of the proof more rigorous.

\bibliography{refer}

\end{document}